\documentclass[12pt]{amsart}
\usepackage{amsmath,amstext,amsfonts,amsthm,amssymb}
\usepackage{eucal}
\usepackage{graphicx}
\renewcommand{\subsubsection}[1]{\addtocounter{subsubsection}{1}
{\ \\[3pt]\bf \thesubsubsection. \  #1} }

\swapnumbers

{  \theoremstyle{definition}

}
\newcommand{\ad}{\operatorname{ad}}

\newcommand{\cl}{{\operatorname{cl}}}

\newcommand{\Ker}{\operatorname{Ker}}

\newcommand{\hra}{\hookrightarrow}

\newcommand{\lra}{\longrightarrow}

\newcommand{\bea}{\begin{eqnarray*}}
\newcommand{\eea}{\end{eqnarray*}}
\newcommand{\bean}{\begin{eqnarray}}
\newcommand{\eean}{\end{eqnarray}}

\newcommand{\ta}{\tilde a}
\newcommand{\tb}{\tilde b}
\newcommand{\tf}{\tilde f}
\newcommand{\tg}{\tilde g}

\newcommand{\tx}{{\tilde{x}}}

\newcommand{\ty}{\tilde y}

\newcommand{\fg}{\mathfrak g}

\newcommand{\CBV}{\mathcal{BV}}

\newcommand{\CT}{\mathcal{T}}

\newcommand{\BN}{\mathbb{N}}

\newcommand{\BZ}{\mathbb{Z}}

\newcommand{\nc}{\newcommand}

\nc{\Id}{\text{Id}}
\nc{\la}{\lambda}

\begin{document}


\centerline{\bf DE RHAM COMPLEX OF A GERSTENHABER ALGEBRA}

\bigskip\bigskip

\centerline{Vadim Schechtman}

\vskip 1cm

\hspace{5cm} {\it \`A Serezha  et Micha, en t\'emoignage d'amiti\'e} 



\vskip 1cm

\centerline{\bf Introduction}

\bigskip\bigskip

The aim of this note is to present some algebraic constructions related to 
Gerstenhaber algebras. 

Let $X$ be a smooth manifold\footnote{we can work in $C^\infty$, analytic, 
or algebraic category}, $\CT_X$ the sheaf of vector fields over $X$. The graded algebra 
of polyvector fields $\Lambda^\bullet\CT_X$ carries a structure of a sheaf  
of {\it Gerstenhaber} (aka {\it odd Poisson)} algebras. This means by definition that $\Lambda^\bullet\CT_X$ carries a Lie bracket of degree $-1$  
(the Schouten - Nijenhuis bracket) such that the operators $\ad_x$ are derivations of the exterior product (with appropriate signs), the precise definitions are recalled below, see 1.1. 

Let us call a {\it Batalin - Vilkovisky (BV)} structure on $\CT_X$ an operator 
$\Delta:\  \Lambda^\bullet\CT_X\lra \Lambda^\bullet\CT_X[-1]$
such that 
$$
\Delta(xy) - \Delta(x)y - (-1)^\tx x\Delta(y) = (-1)^\tx [x,y]
\eqno{(BV1)}
$$
(where we write $\tx = i$ for $x$ being a local section of $\Lambda^i\CT_X$),
and
$$
\Delta\circ \Delta = 0.
\eqno{(BV2)}
$$
According to Koszul [K] (cf. also [S1]),  there is a canonical bijection between the set 
$BV(\CT_X)$ of BV structures 
on $\CT_X$ and the set of integrable connections on the canonical bundle 
$\omega_X$.  It follows that  BV structures exist locally, and form a sheaf 
$\CBV_X$ which is a $\Omega^{1,\cl}_X$-torsor\footnote{here ''cl'' means ''closed''};  
its class is equal to $c_1(\CT_X)\in H^1(X;  \Omega^{1,\cl}_X)$, 
cf. [GMS], \S 11. 

In the present note we generalize the above considerations, replacing  
the Schouten - Nijenhuis algebra $\Lambda^\bullet\CT_X$ 
by an arbitrary Gerstenhaber algebra. 

Let $G^\bullet$ be a  Gerstenhaber algebra. We associate to $G^\bullet$ certain complex $\Omega^\bullet(G^\bullet)$ 
which we call the de Rham complex of $G^\bullet$. 
For example  $\Omega^\bullet( \Lambda^\bullet\CT_X)$ is the usual de 
Rham complex of $X$. In fact, $\Omega^\bullet(G^\bullet)$ appears as  
a subcomplex of a bigger complex $\Omega^\bullet(G^\bullet)^\sim$ 
which might be not without interest in itself.

Let us call a {\bf quasi-BV structure} on $G^\bullet$ an operator 
$\Delta:\  G^\bullet \lra G^\bullet[-1]$ 
which satisfies $(BV1)$ and 
$$
\Delta^2(xy) =\Delta^2(x)y + x\Delta^2(y)
\eqno{(qBV2)}
$$
where $\Delta^2 := \Delta\circ \Delta$. 
It is clear that a BV structure is a quasi-BV structure, and 
if $G^\bullet$ is concentrated in nonnegative degrees and 
generated as an algebra by $G^0$ and $G^1$, as for example  $\Lambda^\bullet\CT_X$, then $(qBV2)$ is equivalent to $(BV2)$. We show that 
quasi-BV structures on $G^\bullet$ form a $\Omega^{1,\cl}(G^\bullet)$-torsor; 
it is a generalization of the described above classification of BV structures. 

The idea of the construction is inspired by [S2] (cf. also [S3]). Since $G^\bullet$ 
is an associative algebra, we can form its Hochschild complex; on the other hand, 
since $G^\bullet[1]$ is a Lie algebra, we can form its Chevalley complex; 
a combination of these two complexes is a bicomplex, and $\Omega^\bullet(G^\bullet)^\sim$ appears as the kernel of  the Hochschild 
differential. 

{\bf Remark.} Formula (BV1) has a curious even counterpart introduced by probabilists, cf. [BE]:
$$
L(fg) - L(f)g - fL(g) = \Gamma(f,g)
$$
where $L$ is a differential operator of the second order, for example the  Laplacian in a Riemannian manifold; the right hand side is called the
{\it  carr\'e du champ} operator and is much used in the theory of Markov 
diffusion processes, cf. [BGL].  Note that the BV operator $\Delta$ is a differential operator of the second order 
with respect to the multiplication (since the Lie - Gerstenhaber bracket is a differential operator of the first order). 

The main part of this note has been worked out back in 2005, after a question asked by  Dmitry Tamarkin in Copenhagen; I am very grateful to him, and also 
to V.Ginzburg, V.Hinich, R.Nest, and B.Tsygan for useful 
discussions. I am especially grateful to the referee for the attentive reading of the text and several 
corrections. 

\bigskip\bigskip

\centerline{\bf \S 1. Gerstenhaber and Batalin - Vilkovisky algebras}

\bigskip\bigskip

In this Section we recall some basic definitions and known results. 

We fix a base field $k$ of characteristic $0$; all our (Lie) algebras, 
modules, tensor products, etc. will be over $k$. 

{\bf 1.1.} 
Recall that a {\it Gerstenhaber algebra} is a $\BZ$-graded commutative algebra 
$G^\bullet$ equipped with a bracket
$$
[ , ]: G^i\otimes G^j\lra G^{i+j-1}
$$
which makes the shifted graded module $G^\bullet[1]$ a graded Lie algebra, that is, it satisfies  
$$
[x,y] = - (-1)^{(\tx-1)(\ty-1)} [y,x]
\eqno{(G1)}
$$
and
$$
[x,[y,z]] = [[x,y],z] + (-1)^{(\tx-1)(\ty-1)} [y,[x,z]]
\eqno{(G2)}
$$
Here we set $\tx = i$ for $x\in G^i$. We denote $G^\bullet[1]$ with 
this Lie bracket by $G^{Lie}$.  
The multiplication and the bracket should be compatible in the following sense: 
$$
[x, yz] = [x,y]z + (-1)^{\ty(\tx-1)} y[x,z]
\eqno{(G3)}
$$
This axiom means that the multiplication 
$G^\bullet\otimes G^\bullet\lra G^\bullet$ is a map of $G^{Lie}$-modules, 
where we consider $G^\bullet$ as a $G^{Lie}$-module by means of 
the adjoint representation, cf. [BD], 1.4.18. 

{\bf 1.2. Example: the Schouten - Nijenhuis algebra.} Let $G^\bullet$ be 
an $\BN$-graded Gerstenhaber algebra, i.e. $G^i = 0$ for $i < 0$. 
Then $A := G^0$ is a commutative algebra,  $T:= G^1$ 
is a Lie algebra and a left $A$-module, and $T$ acts, as a Lie algebra, on $A$ by 
$\tau(a) := [\tau, a],\ a\in A, \tau \in T$; the operators 
$a\mapsto \tau(a)$ are derivations of the ring $A$.  Moreover, the compatibilites 
$$
(a\tau)(b) = a(\tau(b))
\eqno{(LA1)}
$$
$$
[\tau, a\nu] = a[\tau,\nu] + \tau(a)\nu,\ 
\eqno{(LA2)}
$$
$a, b\in A, \tau, \nu\in T$ hold true. 
All this means that $T$ is a {\it Lie algebroid} over $A$. Thus, we've got a 
functor
$$
(Gerstenhaber\ algebras) \lra (Lie\ algebroids).
$$
This functor admits a left adjoint. Namely, if $(A,T)$ is a Lie algebroid, there exists a unique Lie bracket, called the {\it Schouten - Nijenhuis bracket}, cf. [K],  on the exterior algebra $\Lambda^\bullet_A(T)$ 
which coincides with the given bracket on $T$ and makes $\Lambda^\bullet_A(T)$ a Gerstenhaber algebra. 

This is an odd analogue of the classical Poisson structure  on the symmetric algebra $Sym_A^\bullet T$ (whose quantization is the algebra of differential 
operators). 

Let $\Delta$ be a BV structure on $\Lambda^\bullet_A(T)$; consider 
its degree $1$ component: 
$$
c: T\lra A.
\eqno{(1.2.1)}
$$
This operator (sometimes called the {\it divergence}) satisfies the following properties: 
$$
c(a\tau) = ac(\tau) + \tau(a)
\eqno{(Div1)}
$$
$$
c([\tau, \nu]) = \tau(c(\nu)) - \nu (c(\tau)),
\eqno{(Div2)}
$$
$a\in A, \tau, \nu\in T$, which are particular cases of $(BV1)$ and $(BV2)$ 
respectively. Conversely, given an operator  (1.2.1) satisfying $(Div1)$ and 
$(Div2)$, there exists a unique BV structure on $\Lambda^\bullet_A(T)$ 
with $\Delta^1 = c$, cf. [K], [S1], and [GMS], \S 11. (In the last paper 
BV structures on Lie algebroids are called ''Calabi - Yau structures''.)

By that reason we will call an operator (1.2.1) satisfying (Div1) and (Div2) 
a BV structure on the Lie algebroid $T$. 

Set $\Omega^1(T) = Hom_A(T,A)$, and
$$
\Omega^{1,\cl}(T) := \{y\in \Omega^1(T)|\ y([\tau,\nu]) -\tau(y(\nu)) - \nu(y(\tau)) = 0\}.
\eqno{(1.2.2)}
$$ 
The set $BV(T)$ of BV structures on $T$ is canonically an  
$\Omega^{1,\cl}(T)$-torsor, {\it the first Chern class of $T$}. 

\bigskip\bigskip

\centerline{\bf \S 2. Hochschild - Chevalley complex of a Lie algebroid}

\bigskip\bigskip 

We recall here the construction from [S2],  Caput 1. In particular 
we interpret the axioms of a BV structure on $\Lambda^\bullet_AT$ 
as a condition that certain cocycle in a certain ''Hochschild - Chevalley'' bicomplex is a coboundary. 

{\bf 2.1.} If $\fg$ is a Lie algebra and  $M$ a $\fg$-module, 
the  Chevalley complex $C^\cdot_{CH}(\fg,M)$ is defined by  
$C^i_{CH}(\fg,M) = Hom(\Lambda^i\fg,M)$, with the differential  
$$
d_{CH}f(\tau_1,\ldots,\tau_{n+1}) = \sum_{i=1}^{n+1}\ 
(-1)^{i+1}\tau_i f(\tau_1,\ldots, \hat\tau_i,\ldots) + 
$$
$$
+ \sum_{1\leq i < j \leq n+1}\ 
(-1)^{i+j}f([\tau_i,\tau_j],\tau_1,\ldots,\hat\tau_i,\ldots,\hat\tau_j, 
\ldots)
\eqno{(2.1.1)}
$$
If $A$ is a commutative algebra, and   $M, N$  two $A$-modules, let us define 
the bar (Hochschild) complex $C^\cdot_H(M,N)$ by 
$C^i_H(M,N) = Hom(A^{\otimes i}\otimes M,N)$ with a diffierential   
$$
d_Hf(a_1,\ldots,a_{n+1},x) = a_1f(a_2,\ldots,a_{n+1},x) + 
\sum_{i=1}^n\ (-1)^if(a_1,\ldots, a_ia_{i+1},\ldots) + 
$$
$$
+ (-1)^{n+1}f(a_1,\ldots, a_n;a_{n+1}x)
\eqno{(2.1.2)}
$$
So for example  $H^0_{H}(M,N) := H^0C^\cdot_H(M,N) = Hom_A(M,N)$. 

{\bf 2.2.} Let $T$ be a Lie algebroid over $A$.  Let us define a bicomplex  
$C^{\cdot\cdot}_{HCH}(T,A)$, called the  {\it Hochschild - Chevalley bicomplex} of $T$ as follows. By definition, the $0$-th column 
will be the Hochschild complex  
$$
C^{0\cdot}_{HCH}(T,A) = C^{\cdot}_{H}(T,A):\ Hom(T,A) 
\lra Hom(A\otimes T,A) \lra \ldots,
$$
and the $0$-th line will be a truncated and shifted Chevalley complex    
$$
C^{\cdot 0}_{HCH}(T,A) = C^{\cdot + 1}_{CH}(T,A):\ Hom(T,A) \lra 
Hom(\Lambda^2 T,A) \lra \ldots, 
$$
with the differential multiplied by $-1$.

For   $i \geq 0,\ j \geq 1$ we set 
$$
C^{ij}_{HCH}(T,A) := 
Hom(A^{\otimes j}\otimes T\otimes \Lambda^i T,A) 
$$
The horizontal differential is the Chevalley one, after the identification 
$$
Hom(A^{\otimes j}\otimes T\otimes \Lambda^i T,A)
= 
Hom(\Lambda^i T, Hom( A^{\otimes j}\otimes T,A)),
$$
explicitly:
$$
d_{DR}f(a_1,\ldots,a_j;\tau;\tau_1,\ldots ) = \sum_p\ (-1)^{p+1}\{ 
\tau_p f(a_1,\ldots;\tau;\ldots,\hat\tau_p,\ldots) - 
$$
$$ 
- \sum_r\ f(a_1,\ldots,\tau_p(a_r),\ldots;\tau;\ldots,\hat\tau_p,\ldots) 
- f(a_1,\ldots;[\tau_p,\tau];\ldots,\hat\tau_p,\ldots)\} +    
$$
$$
+ \sum_{p<q}\ (-1)^{p+q} 
f(a_1,\ldots;\tau;\ldots,\hat\tau_p,\ldots,\hat\tau_q,\ldots) 
$$
The vertical differentials
$$
d_H:\ C^{i-1,0} = Hom(\Lambda^iT,A) \lra 
Hom(A\otimes T\otimes \Lambda^{i-1}T,A) = C^{i-1,1}
$$
are given by
$$
d_Hf(a;\tau_1,\ldots,\tau_i) = af(\tau_1,\ldots,\tau_i) - 
f(a\tau_1,\ldots,\tau_i)
$$
where we identify $Hom(\Lambda^iT,A)$ with the space of 
skew-symmetric maps $T^{\otimes i}\lra A$. 

Note that
$$
\Ker d_H = \Omega^i(T):= Hom_A(\Lambda^i_AT, A)
$$
(a skew-symmetric map $A$-linear with respect to the first 
argument is $A$-linear with respect to all the other ones). 

One checks by hand that the squares 
$$
\begin{matrix} C^{i1} & \overset{d_{CH}}\lra & C^{i+1,1}\\
d_H\uparrow & & \uparrow d_H\\
C^{i0} & \overset{d_{CH}}\lra & C^{i+1,0}
\end{matrix}
$$
are commutative. 

The vertical differential $C^{i-1,j}\lra C^{i-1,j+1}$ for $j\geq 1$ 
is the Hochschild's one, with respect to the first argument:
$$
d_Hf(a_1,\ldots,a_j;\tau_1;\tau_2,\ldots,\tau_i) = 
a_1 f(a_2,\ldots,a_j;\tau_1;\tau_2,\ldots,\tau_i) + 
$$
$$
+ \sum_{p=1}^{j-1}\ (-1)^p 
f(a_1,\ldots,a_pa_{p+1},\ldots,a_j;\tau_1;\tau_2,\ldots,\tau_i) +
$$
$$
+ (-1)^j f(a_1,\ldots,a_{j-1};a_j\tau_1;\tau_2,\ldots,\tau_i) 
$$ 
The commutation     
$d_Hd_{CH} = d_{CH}d_H$ for upper squares follows from  

{\bf 2.3.} {\it Lemma.} The  Hochschild differentials 
$$
d_H:\ Hom(A^{\otimes n}\otimes T,A) \lra 
Hom(A^{\otimes (n+1)}\otimes T,A)
$$
are morhisms of  $T^{Lie}$-modules.  $\square$
$$
\begin{matrix}
\uparrow& & \uparrow& &\uparrow& \\
Hom(A^{\otimes 2}\otimes T, A)&\lra &
Hom(A^{\otimes 2}\otimes T\otimes T, A)&\lra & 
Hom(A^{\otimes 2}\otimes T\otimes \Lambda^2T, A)&\lra\\
\uparrow& & \uparrow& &\uparrow& \\
Hom(A\otimes T, A)&\lra &
Hom(A\otimes T\otimes T, A)&\lra & 
Hom(A\otimes T\otimes \Lambda^2T, A)&\lra\\
\uparrow& & \uparrow& &\uparrow& \\
Hom(T, A)&\lra &
Hom(\Lambda^2T, A)&\lra & 
Hom(\Lambda^3T, A)&\lra\\ 
\end{matrix}
$$ 

{\bf 2.4.} So one has defined a double complex $C^{\cdot\cdot}_{HCH}(T,A)$ (drawn above), 
whence the associated simple complex  $C^\cdot_{HCH}(T,A)$, with the 
differential defined by the usual formula 
$$
d_{HCH}(x^{ij}) = d_{CH}(x^{ij}) + (-1)^i d_H(x^{ij}) 
$$
As we have remarked above,   
$$
H^0_H(C^{i\cdot}_{CHC}(T,A)) = \Omega^{i+1}(T)
$$
whence a canonical emebedding of the (truncated and) shifted  de Rham complex 
$$
\Omega^{\cdot + 1}(T) \hra C^\cdot_{HCH}(T,A).
$$
In particular
$$
H^0 (C^\cdot_{HCH}(T,A)) = \Omega^{1,\cl}(T) := \text{Ker}(d_{DR}:\ 
\Omega^1(T) \lra \Omega^2(T))
$$
cf. (1.2.2). 

{\bf 2.5.} Consider a canonical element  
$$
e \in C^{01}_{HCH}(T,A) = Hom(A\otimes T,A),\ e(a,\tau) = \tau(a).
$$
One has 
$$
d_He(a,b;\tau) = a\tau(b) - \tau(ab) + b\tau(a) = 0.
$$
Similarly
$$
d_{CH}e(\tau,a,\tau') = Lie_{\tau'}e(a,\tau) = 
\tau'\tau(a) -  \tau\tau'(a) - [\tau',\tau](a) = 0.
$$
It follows that the element
$$
\epsilon = (-e,0) \in C^{01}_{HCH}(T,A) \oplus C^{10}_{HCH}(T,A) = 
C^{1}_{HCH}(T,A)
$$
is a $1$-cocycle in the total complex   $C^{\cdot}_{HCH}(T,A)$. 

When $\epsilon$ is a coboundary? Equation 
$$
d_{HCH} c = \epsilon, \ c\in C^{0}_{HCH}(T,A) = Hom(T,A)
$$
is equivalent to two equations: 
$$
d_Hc = - e, 
\eqno{(2.5.1)}
$$
i.e. 
$$
ac(\tau) - c(a\tau) = - \tau(a)
$$
which is  the first axiom 1.2 (Div1), and 
$$
d_{CH} c = 0, 
\eqno{(2.5.2)}
$$
i.e. 
$$
\tau c(\tau') - \tau' c(\tau) + c([\tau,\tau']) = 0 
$$
which is 1.2 (Div2). It follows that a BV structure on $T$ is the same as an element  
 $c\in C^0_{HCH}(T,A)$ such that $d_{HCH} c = \epsilon$.   

Consequently, the set   
$$
CY(T) = \{ c\in C^0_{HCH}(T,A) | d_{HCH} c = \epsilon\} 
$$
is naturally a torsor under  
$$
H^0 C^\cdot_{HCH}(T,A) = \Omega^{1,\cl}(T),
$$
as we have seen already previously.

\bigskip\bigskip 

\centerline{\bf \S 3.  
Hochschild complex of a Gerstenhaber algebra } 


\bigskip\bigskip

{\bf 3.1.} Let  $A^\cdot = \oplus_{i\in \BZ} A^i$ be a graded $k$-module. 
For $a \in A^i$ we set
$$
\ta = |a| = i 
$$
If $B^\cdot$ is another graded module, we set 
$$
(a\otimes b)^\sim = \ta +  \tb
$$
for $a \otimes b \in A^\cdot \otimes B^\cdot$. If $f:\ A^\cdot \lra B^\cdot$ is a map, we say that 
 $\tf = i$ if
$$
f(a)^\sim = \ta + i
$$ 
Thus for  $f \in Hom(A^{\cdot\otimes n},A^\cdot)$, we say that $\tf = i$ if
$$
f(a_1,\ldots,a_n)^\sim = \sum_{j=1}^n\ \ta_j + i
$$
We denote 
$$
Hom^i(A^\cdot,B^\cdot) = \{f\in Hom(A^\cdot,B^\cdot) | \tf = i\}
$$

{\bf 3.2.} Let  $A^\cdot$ be an associative graded $k$-algebra,  
$M^\cdot$ a graded $A^\cdot$-bimodule. One defines the Hochschild complex 
$C^\cdot_H(A^\cdot,M^\cdot)$ by 
$$
C_H^n(A^\cdot,M^\cdot) = Hom(A^{\cdot\otimes n},M^\cdot)
$$
For $f \in C_H^n(A^\cdot,M^\cdot)$ set 
$$
|f| = \tf + n
$$
(sic!). The differential is defined by the formula  
(cf. [TT], 2.3):
$$
d_H f(a_1,\ldots,a_{n+1}) = (-1)^{|a_1||f| + |f| + 1} 
a_1 f(a_2,\ldots,a_{n+1}) + 
$$
$$
+ \sum_{i=1}^n\ (-1)^{|f| + 1 + \sum_{p = 1}^i (|a_p|+1)} 
f(a_1,\ldots,a_ia_{i+1},\ldots, a_{n+1}) + 
$$
$$
+ (-1)^{|f| + \sum_{p = 1}^n (|a_p|+1)} 
f(a_1,\ldots, a_n)a_{n+1}
\eqno{(3.2.1)}
$$
One checks that $d_H^2 = 0$. 

We will need the case $M^\cdot = A^\cdot$. 

{\bf 3.3.} For example,  
$$
d_Hf(a,b) = (-1)^{|a||f| + |f| + 1}a f(b) + 
$$
$$
+ (-1)^{|f| + |a|} f(ab) + (-1)^{|f| + |a| + 1}f(a)b
\eqno{(3.3.1)}
$$ 
Similarly    
$$
d_Hf(a,b,c) = (-1)^{|a||f| + |f| + 1}a f(b,c) + 
(-1)^{|f| + |a|} f(ab,c) +
$$
$$
+ (-1)^{|f| + |a| + |b| + 1} f(a,bc) 
+ (-1)^{|f| + |a| + |b|} f(a,b)c.
\eqno{(3.3.2)}
$$
In particular, the kernel
$$
\text{Ker} (d_H^1:\ C^1_H(A^\cdot,M^\cdot) \lra C^2_H(A^\cdot,M^\cdot) = 
$$
$$
= \{ f:\ A^\cdot \lra M^\cdot |\ - (-1)^{\tf\ta} af(b) + f(ab) - 
f(a)b = 0\} = 
$$
$$
= \text{Der} (A^\cdot, M^\cdot)
\eqno{(3.3.3)}
$$
is identified with the module of derivations of  $A^\cdot$ with values in 
$M^\cdot$. 

{\bf 3.4.} Let $G^\cdot$ be a Gerstenhaber algebra, and consider its  Hochschild $C_H^\cdot(G^\cdot,G^\cdot)$ 
complex. Let us consider the Lie bracket
$[ , ]:\ G^\cdot\otimes G^\cdot \lra G^\cdot$ as a Hochschild cochain 
 $[ , ]\in C^2_H(G^\cdot,G^\cdot)$. We have 
$\deg\ [ , ] = -1$, therefore 
$$
| [ , ] | = 1.
$$

{\bf 3.5.} {\bf Lemma.} $d_H [ , ] = 0$. 

{\bf Proof.} We will use the Poisson identity
$$
[x,yz] = [x,y]z + (-1)^{|y|(|x| + 1)}y[x,z] 
\eqno{(3.5.1)}
$$
and its  companion 
$$
[xy,z] = x[y,z] + (-1)^{|y|(|z| + 1)} [x,z]y
\eqno{(3.5.2)}
$$
which is a consequence of (3.5.1) and of the  commutativity. 

Using  (3.3.2), 
$$
d_H[ , ](x,y,z) = (-1)^{|x|}x[y,z] + (-1)^{|x| + 1}[xy,z] + 
$$
$$
+ (-1)^{|x| + |y|}[x,yz] + (-1)^{|x| + |y| + 1}[x,y]z = 
$$
$$
= (-1)^{|x|}x[y,z] + 
(-1)^{|x| + 1}\{x[y,z] + (-1)^{|y|(|z| + 1)} [x,z]y\} + 
$$
$$
+ (-1)^{|x| + |y|}\{[x,y]z + (-1)^{|y|(|x| + 1)}y[x,z]\} + 
(-1)^{|x| + |y| + 1}[x,y]z = 0,
$$ 
{\it QED}. 

{\bf 3.6.} {\bf Remark.} Maybe the identity  $d_H[, ] = 0$, i.e. 
$$
x[y,z] - [xy,z] + 
(-1)^{|y|}[x,yz] - (-1)^{|y|}[x,y]z = 0
$$
can serve as a replacement of the Poisson identity 
for non-commutative Poisson algebras.

{\bf 3.7.} For
$$
\Delta \in Hom^{-1}(G^\cdot,G^\cdot) \subset C^1(G^\cdot,G^\cdot) 
$$
we will have $|\Delta| = 0$, and  
$$
d_H\Delta(x,y) = - x\Delta(y) + (-1)^{|x|}\Delta(xy) + 
(-1)^{|x| + 1}\Delta(x)y, 
$$
cf. (3.3.1).  Therefore the equation 
$$
d_H\Delta = [ , ]
\eqno{(3.7.1)}
$$
is equivalent to 
$$
\Delta(xy) - \Delta(x)y - (-1)^{|x|} x\Delta(y) = (-1)^{|x|}[x,y], 
$$
that is, to  the axiom (BV1), cf. 1.3 (Koszul would say that  $\Delta$ {\it generates} the bracket $[ , ]$, cf. [K], p. 262).  

{\bf 3.8.} Set  $A = G^0,\ T = G^1$.  One has a canonical morphism of Gerstenhaber algebras 
$$
\Lambda^\cdot_A(T) \hra G^\cdot.
\eqno{(3.8.1)}
$$
Evidently,  $Hom(T,A) \subset Hom^{-1}(G^\cdot,G^\cdot)$. 
One has a canonical projection
$$
\pi:\ Hom(G^\cdot,G^\cdot) \lra Hom(T,A)
$$
which may be extended to a canonical projection of the
Hochschild complexes
$$
\pi:\ C_H^{\cdot - 1}(G^\cdot,G^\cdot) \lra C_H^\cdot(T,A).
$$
Note that the last term of the Hochschild differential  
 (3.2.1) disappears in $C_H^\cdot(T,A)$.  

One has
$$
\pi([ , ]) = - e, 
$$
since $[a,\tau] = - \tau(a)$, and an operator  $\Delta$ satisfying  
(3.7.1) will be transformed into an operator  $c\in Hom(T,A)$ satisfying  $(Div1)$, $d_H c = - e$.  

\bigskip\bigskip

\centerline{\bf \S 4. The de Rham complex of a Gerstenhaber algebra} 
\bigskip\bigskip 

{\bf 4.1.} Let $\fg^\cdot$ be a graded Lie algebra and 
$M^\cdot$ a graded $\fg^\cdot$-module. The  Chevalley complex 
$C^\cdot_{CH}(\fg^\cdot,M^\cdot)$ is  defined by 
$C^n(\fg^\cdot,M^\cdot) = Hom(\Lambda^n\fg^\cdot,M^\cdot)$,   
with a differential
$$
d_{CH}f(x_1,\ldots,x_{n+1}) = 
$$
$$
= \sum_{i=1}^{n+1}\ 
(-1)^{i + 1 + \tx_i(\tf + \sum_{p=1}^{i-1} \tx_p)} 
x_if(x_1,\ldots,\hat x_i,\ldots) + 
$$
$$
+ \sum_{1 \leq i < j \leq n+1}\ 
(-1)^{i+j + \tx_i\sum_{p=1}^{i-1} \tx_p + 
\tx_j\sum_{1 \leq q \leq j-1, q\neq i} \tx_q}\times 
$$
$$ 
\times f([x_i,x_j],x_1,\ldots,\hat x_i, 
\ldots, \hat x_j,\ldots)
$$
If $N^\cdot$ is another graded  $\fg^\cdot$-module,  $\fg^\cdot$ 
acts on  $Hom(M^\cdot,N^\cdot)$ by the usual formula
$$
(\tau\phi)(x) = \tau (\phi(x)) - (-1)^{\deg\ \phi\tx}.
\phi(\tau x)
$$

{\bf 4.2.} Let $G^\cdot$ be a Gerstenhaber algebra. We can apply the previous definition 
to the graded Lie algebra 
$G^{Lie\cdot} := G^\cdot[1]$ which acts on itself by the adjoint representation, and get the complex 
$$
C^\cdot_{CH}(G^\cdot,G^\cdot) := 
C^\cdot_{CH}(G^{Lie\cdot},G^{Lie\cdot}).
$$
More explicitly, 
$$
C^n_{CH}(G^\cdot,G^\cdot) = Hom(\Lambda^n(G^{\cdot}[1]),G^{\cdot}[1])
\subset
$$
$$  
\subset Hom(G^{\cdot}[1])^{\otimes n},G^{\cdot}[1])
= Hom(G^{\cdot\otimes n},G^{\cdot})
$$
For a homogenenous element
 $f\in C^n_{CH}(G^\cdot,G^\cdot)$, let us denote by 
 $\tf$ its degree as an element of  
$Hom(G^{\cdot\otimes n},G^{\cdot})$. 

For an element $x\in G^i$, its degree as an element of 
$G^{\cdot}[1]$ is $i - 1$. 

Therefore, given a map 
$f:\ (G^{\cdot})^{\otimes n} \lra G^{\cdot}$ 
of degree $\tf$, its degree as an element of 
$Hom(G^{\cdot}[1])^{\otimes n},G^{\cdot}[1])$ is equal to 
$\tf + n - 1$.  

One identifies 
 $C^n_{CH}(G^\cdot,G^\cdot)$ with the space of polylinear  functions
 $f:\ (G^\cdot)^n \lra G^\cdot$ which are {\it alternating} in the following sense: 
$$
f(x_1, \ldots, x_i,x_{i+1},\ldots,x_n) = 
$$
$$
= - (-1)^{(\tx_i - 1)(\tx_{i+1} - 1)} 
f(x_1, \ldots, x_{i+1},x_{i},\ldots,x_n)
\eqno{(4.2.1)}
$$ 
The differential is acting as follows: for
$f \in C^n_{CH}(G^\cdot,G^\cdot)$,  
$$
d_{CH}f(x_1,\ldots,x_{n+1}) = \sum_{i=1}^{n+1}\ 
[x_i,f(x_1,\ldots,\hat x_i,\ldots)] \times 
$$
$$
\times (-1)^{i+1+(\tx_i - 1)(\tf + n - 1 + 
\sum_{p=1}^{i-1} (\tx_p - 1))} + 
$$
$$
+ \sum_{1\leq i < j \leq n+1}\ f([x_i,x_j],x_1,\ldots,\hat x_i,\ldots, 
\hat x_j,\ldots)\times 
$$   
$$
\times (-1)^{i+j + (\tx_i - 1)\sum_{p=1}^{i-1} (\tx_p - 1)}
\times 
$$
$$
\times (-1)^{(\tx_j - 1)\sum_{1 \leq q \leq j - 1,\ q\neq i}\ 
(\tx_q - 1)}.
\eqno{(4.2.2)}
$$
For example, for  $g\in G^\cdot$ we have  
$$
d_{CH}g(x) = (-1)^{(\tx - 1)(\tg - 1)}[x,g] = [g,x],
\eqno{(4.2.3)}
$$
and for $f\in Hom(G^\cdot,G^\cdot)$ 
$$
d_{CH}f(x,y) = [x,f(y)](-1)^{(\tx - 1)\tf} - 
$$
$$ 
- [y,f(x)](-1)^{(\ty - 1)(\tf + \tx - 1)} 
- f([x,y]) = 
$$
$$
= (-1)^{(\tx - 1)\tf}[x,f(y)] + [f(x),y] - f([x,y]).
\eqno{(4.2.4)} 
$$
Note that 
 $(d_{CH}f)^\sim = \tf - 1$. 

{\bf 4.3.} Now we will introduce a double ''Hochschild - Chevalley'' complex 
$\{ C^{ij}_{HCH}(G^\cdot,G^\cdot)\}$, with $i \geq 0,\ j = 0, 1$.  
So our double complex will have only two lines. 

The $0$-th line will be a shifted Chevalley complex:  
$$
C^{\cdot 0}_{HCH}(G^\cdot,G^\cdot) := 
C^{\cdot + 1}_{CH}(G^\cdot,G^\cdot).
$$
For example, at the corner we will have 
$$
C^{00}_{HCH}(G^\cdot,G^\cdot) = Hom(G^\cdot,G^\cdot). 
$$

{\bf 4.4.} The first line of the Hochschild - Chevalley double complex coincides by definition with the  
 Chevalley complex of $G^{Lie\cdot}$ with coefficients in  
$C^{2}_H(G^\cdot,G^\cdot)$: 
$$
C^{\cdot 1}_{HCH}(G^\cdot,G^\cdot) := 
C^\cdot_{CH}(G^{Lie},C^{2}_H(G^\cdot,G^\cdot)),\ n \geq 1. 
$$
More exactly, one identifies  $C^{2}_H(G^\cdot,G^\cdot)$ with  
$$
Hom(G^{\cdot\otimes 2},G^\cdot) = 
Hom((G^{\cdot}[1])^{\otimes 2},G^\cdot[1]), 
$$
where the right hand side carries an evident structure of a 
 $G^{Lie}$-module. 

Note that one can identify
$$
C^{m 1}_{HCH}(G^\cdot,G^\cdot) = Hom(\Lambda^m(G^{\cdot}[1])\otimes
(G^{\cdot}[1])^{\otimes 2},G^\cdot[1]) 
$$
which in turn is a spaceof functions 
of $m+2$ arguments
$$
f(x_1,\ldots,x_m;x_{m+1},x_{m+2}):\ (G^\cdot)^{m+2} \lra G^\cdot
$$
which are alternating, in the sense of 
 (4.2.1), with respect to the first  $m$ 
arguments. 

{\bf 4.5.} Explicitly, the Chevalley differential in the first line  
$$
d_{CH}:\ C^{n-1,1}_{HCH}(G^\cdot,G^\cdot) = 
Hom(\Lambda^{n-1}(G^{\cdot}[1])\otimes
(G^{\cdot}[1])^{\otimes 2},G^\cdot[1]) \lra
$$
$$
\lra Hom(\Lambda^{n}(G^{\cdot}[1])\otimes
(G^{\cdot}[1])^{\otimes 2},G^\cdot[1]) = 
C^{n1}_{HCH}(G^\cdot,G^\cdot)
$$
acts by 
$$
d_{CH}f(x_1,\ldots,x_n;x_{n+1},x_{n+2}) = 
$$     
$$
= \sum_{i=1}^n\ [x_i,f(x_1,\ldots,\hat x_i,\ldots)]\cdot 
(-1)^{i+1 + (\tx_i - 1)[\tf + n + \sum_{p=1}^{i-1} (\tx_p - 1)]} - 
$$
$$
- f(x_1,\ldots,\hat x_i,\ldots; [x_i,x_{n+1}],x_{n+2})\cdot 
(-1)^{i + 1 +(\tx_i - 1)\sum_{p=i+1}^n (\tx_p - 1)} - 
$$
$$
- f(x_1,\ldots,\hat x_i,\ldots; x_{n+1}, [x_i,x_{n+2}])\cdot 
(-1)^{i + 1 +(\tx_i - 1)\sum_{p=i+1}^{n+1} (\tx_p - 1)} + 
$$
$$
+ \sum_{1\leq i < j \leq n}\ f([x_i,x_j],x_1,\ldots,\hat x_i,\ldots, 
\hat x_j,\ldots;x_{n+1},x_{n+2})\times 
$$
$$
\times (-1)^{i+j + (\tx_i - 1)\sum_{p=1}^{i-1} (\tx_p - 1)}
\times 
$$
$$
\times (-1)^{(\tx_j - 1)\sum_{1 \leq q \leq j - 1,\ q\neq i}\ 
(\tx_q - 1)}
\eqno{(4.5.1)}
$$
For example, for $n=1$: 
$$
d_{CH}f(x;y,z) = [x,f(y,z)](-1)^{(\tx - 1)(\tf + 1)} - 
$$
$$
- f([x,y],z) - f(y,[x,z])(-1)^{(\tx - 1)(\ty - 1)}. 
\eqno{(4.5.2)} 
$$

\bigskip

(a) {\it The corner square}

\bigskip

{\bf 4.6.} Recall that 
$$
d_H = d_H^{00}:\  Hom(G^\cdot[1],G^\cdot[1]) \lra 
Hom((G^\cdot[1])^{\otimes 2},G^\cdot[1])
$$
acts as 
$$
d_Hf(x,y) = (-1)^{\tx(\tf + 1) + \tf}xf(y) + 
$$
$$
+ (-1)^{\tf + \tx + 1}f(xy) 
+ (-1)^{\tf + \tx}f(x)y, 
$$
cf. (3.3.1) (we have $|f| = \tf + 1$). 

One defines a vertical ''Hochschild'' differential 
$$
d_H = d_H^{10}:\ C^{10}_{HCH}(G^\cdot,G^\cdot) = 
Hom(\Lambda^{2}(G^\cdot[1]),G^\cdot[1]) \lra 
$$
$$
\lra Hom((G^\cdot[1])^{\otimes 3},G^\cdot[1]) = 
C^{11}_{HCH}(G^\cdot,G^\cdot)
$$
by the formula 
$$
d_Hf(x,y,z) = (-1)^{\ty(\tf + \tx + 1) + \tf + \tx}
yf(x,z) + 
$$
$$
+ (-1)^{\tf + \tx + \ty + 1}f(x,yz) 
+ (-1)^{\tf + \tx + \ty}f(x,y)z 
$$

{\bf 4.7.} {\bf  Lemma.} {\it We have 
$$
d_{CH}^{01}d_H^{00} = d_H^{10}d_{CH}^{00}.
$$}

This is proven by a direct (although long) calculation.

\bigskip 

(b) {\it The vertical differential: general case }

\bigskip

{\bf 4.8.} One defines the vertical ''Hochschild'' differential 

$$
d_H^{n0}:\ C^{n0}_{HCH}(G^\cdot,G^\cdot) = 
Hom(\Lambda^{n+1}(G^\cdot[1]),G^\cdot[1]) \lra 
$$
$$
\lra Hom(\Lambda^{n}(G^\cdot[1])\otimes (G^\cdot[1])^{\otimes 2} ,G^\cdot[1]) 
= C^{n1}_{HCH}(G^\cdot,G^\cdot)
$$
by the formula 
$$
d_H^{n0}f(x_1,\ldots,x_n;x_{n+1},x_{n+2}) = 
(-1)^{\tx_{n+1}(\tf + \sum_1^n\ \tx_p + 1) + 
\tf + \sum_1^n\ \tx_p} \times 
$$
$$
\times x_{n+1}f(x_1,\ldots,x_n,x_{n+2}) + 
$$
$$  
+ (-1)^{\tf + \sum_1^{n+1}\ \tx_p + 1}
f(x_1,\ldots,x_n,x_{n+1}x_{n+2}) + 
$$
$$
+ (-1)^{\tf + \sum_1^{n+1}\ \tx_p}
f(x_1,\ldots,x_n,x_{n+1})x_{n+2}, 
$$
Note that
$$
(d_{CH}f)^\sim = \tf - 1\ \text{et\ }(d_{H}f)^\sim = \tf. 
$$

{\bf 4.9.} {\bf Lemma.} {\it For all  $n \geq 0$, on a $d_H^{n0}d_{CH}^{n-1,0} = 
d^{n-1,1}_{CH}d^{n-1,0}_H$.}

The proof is also a direct  calculation.   

\bigskip

(c) {\it De Rham complex} 

\bigskip

{\bf 4.13.} We define  spaces 
$$
\Omega^n(G\cdot)^\sim := \text{Ker}(d_H:\ 
C^{n-1,0}_{HCH}(G^\cdot) \lra C^{n-1,1}_{HCH}(G^\cdot)),\ n\geq 1; 
$$
$\Omega^0(G\cdot)^\sim = G^\cdot$. 

Explicitly, the space
$\Omega^n(G\cdot)^\sim$ is identified with the space of polylinear maps 
$$
f:\ G^{\cdot n} \lra G^\cdot
$$
which are:

(a) alternating, i.e. 
$$
f(x_1,\ldots,x_n) = (-1)^{(\tx_i - 1)(\tx_{i+1} - 1)} 
f(x_1,\ldots,x_{i+1},x_i,\ldots,x_n)
$$

(b) satisfy the identity
$$
f(x_1,\ldots,x_{n-1},yz) = 
$$
$$
= (-1)^{y(\tf + \sum_1^{n-1} \tx_p)} 
yf(x_1,\ldots,x_{n-1},z) + f(x_1,\ldots,x_{n-1},y)z
$$

For example, 
$$
\Omega^1(G^\cdot)^\sim = \text{Der}(G^\cdot,G^\cdot) := 
\{ f \in Hom(G^\cdot,G^\cdot)|\ 
f(xy) = (-1)^{\tx\tf} xf(y) + f(x)y\} 
$$
Due to 4.9, the Chevalley differential induces a differential  
$d_{DR}:\ \Omega^n(G^\cdot) \lra \Omega^{n+1}(G^\cdot)$, $n\geq 1$. 

Moreover, recall the differential
$$
d_{CH}:\ G^\cdot \lra Hom(G^\cdot,G^\cdot),\  
d_{CH}g(x) = [g,x],
$$
 cf. 4.2.3. We have $(d_{CH}g)^\sim = \tg - 1$, and one checks at once
that  $d_{CH}g \in \text{Der}(G^\cdot,G^\cdot)$, therefore
$d_{CH}$ induces a map 
$$
d_{DR}:\ G^\cdot = \Omega^0(G^\cdot)^\sim \lra \Omega^1(G^\cdot)^\sim = 
\text{Der}(G^\cdot,G^\cdot),\ d_{DR}g(x) = [g,x].
$$
This defines a complex $\Omega^\bullet(G^\bullet)^\sim$ which we call 
a {\bf big de Rham complex of $G^\bullet$}. 

Inside it, we define a subcomplex
$\Omega^\cdot(G^\cdot) \subset \Omega^\cdot(G^\cdot)^\sim$, 
{\bf the small de Rham complex} of $G^\cdot$, by 
$$
\Omega^0(G^\cdot) = G^0;\ 
\Omega^n(G^\cdot) = \{f \in \Omega^n(G^\cdot)^\sim | \tf = - n\}
$$

{\bf 4.14. Example.} Suppose that $G^\cdot = \Lambda_A^\cdot(T)$ is the  
Schouten - Nijenhuis algebra of an  $A$-Lie algebroid  $T$. We have
$\Omega^0(\Lambda_A^\cdot(T)) = A$. 

A degree $-1$ derivation  
$\omega \in \Omega^1(\Lambda_A^\cdot(T)) = 
\text{Der}^{-1}(\Lambda_A^\cdot(T), \Lambda_A^\cdot(T))$ is uniquely defined 
by its values on  $T$, $\omega(\tau) \in A$. 
The derivation condition for  $a\in A,\ \tau \in T$ will be  
$\omega(a\tau) = a\omega(\tau)$, i.e. 
$$
\Omega^1(\Lambda_A^\cdot(T)) = Hom_A(T,A) = \Omega^1(T).
$$
Similarly, a biderivation 
 $\Omega\in \Omega^2(\Lambda_A^\cdot(T))$ is uniquely defined by its values   
$\omega(\tau,\tau') \in A,\ 
\tau, \tau' \in T$, and the axiom (b) means that  $\omega$ is  $A$-linear with respect to  
 $\tau'$, and therefore with respect to 
 $\tau$, since 
$\omega$ is alternating, whence  
$$
\Omega^2(\Lambda_A^\cdot(T)) = Hom_A(\Lambda^2T,A) = \Omega^2(T).
$$
More generally,   $\Omega^n(\Lambda_A^\cdot(T))$ may be identified with 
$\Omega^n(T):= Hom_A(\Lambda^n T,A)$ for all pour tous $n\geq 0$.  

On the other hand, for $\omega\in \Omega^1(\Lambda_A^\cdot(T))$ and  
$\tau, \tau'\in T$ we have (cf. (4.2.4)):
$$
d_{DR}\omega(\tau,\tau') = d_{CH}\omega(\tau,\tau') = 
[\tau,\omega(\tau')] - [\tau',\omega(\tau)] 
- \omega([\tau,\tau']) = 
$$
$$
= \tau(\omega(\tau')) - \tau'(\omega(\tau) 
- \omega([\tau,\tau'])
$$
More generally, the  differential in  
$\Omega^\cdot(\Lambda_A^\cdot(T))$ is identified with the usual de Rham differential in  $\Omega^\cdot(T)$, 
cf. [GMS], 1.3. 

In particular, the subspace of closed $1$-forms
$$
\Omega^{1,\text{cl}}(\Lambda_A^\cdot(T)) = 
\Omega^{1,\text{cl}}(T).
$$

\bigskip\bigskip

\centerline{\bf \S 5. Batalin - Vilkovisky structures} 

\bigskip\bigskip

{\bf 5.1.} Let $G^\cdot$ be a Gerstenhaber algebra. Consider the bracket   
$$
[,] \in Hom(G^{\cdot\otimes 2},G^\cdot) = C^{01}_{CHC}(G^\cdot,G^\cdot),\ 
[,] = - 1 
$$
We have seen in 3.4 that $d_H[ , ] = 0$. 

On the other hand (cf. (4.5.2)), 
$$
d_{CH}[ , ](x,y,z) = [x,[y,z]] - [[x,y],z] - (-1)^{(\tx - 1)(\ty - 1)} 
[y,[x,z]] = 0, 
$$
i.e. $d_{CH}[ , ] = 0$ as well. 

{\bf 5.2.} Let  
$\Delta \in Hom(G^\cdot,G^\cdot)$ be an element of degree  
$-1$ such that 
$$
d_H\Delta = [ , ], 
$$
i.e. (cf. 3.7) 
$$
[x,y] = (-1)^{\tx}\{\Delta(xy) - \Delta(x)y - (-1)^{\tx} x\Delta(y)\}
\eqno{(5.2.1)} 
$$
Let us compute $d_{CH}\Delta$ (cf. (4.2.4)): 
$$
d_{CH}\Delta(x,y) = (-1)^{\tx - 1}[x,\Delta(y)] 
- (-1)^{(\ty - 1)\tx}[y,\Delta(x)] - \Delta([x,y]) = 
$$ 
$$
= - \biggl\{\Delta(x\Delta(y)) - \Delta(x)\Delta(y) 
- (-1)^{\tx}x\Delta^2(y)\biggr\} - 
$$
$$
- (-1)^{(\ty - 1)\tx + \ty}\biggl\{\Delta(y\Delta(x)) - \Delta(y)\Delta(x) 
- (-1)^{\ty} y\Delta^2(x)\biggr\} - 
$$
$$
- (-1)^{\tx}\biggl\{\Delta^2(xy) - \Delta(\Delta(x)y) 
- (-1)^{\tx}\Delta(x\Delta(y))\biggr\} = 
$$
$$
= (-1)^{\tx + 1}\{\Delta^2(xy) - \Delta^2(x)y - x \Delta^2(y)\}
$$
Note that 
$$
(-1)^{(\ty - 1)\tx}[y,\Delta(x)] = [\Delta(x),y]
$$

Thus  we have:

{\bf 5.3.} {\bf Lemma.} {\it If $\Delta \in Hom(G^\cdot,G^\cdot)$, $\Delta^\sim = 
-1$ v\'erifie $d_H\Delta = 0$ then
$$
d_{CH}\Delta(x,y) = 
(-1)^{\tx + 1}\{\Delta^2(xy) - \Delta^2(x)y - x \Delta^2(y)\}.
$$}
$\square$

{\bf 5.4.} Let us call a {\bf quasi BV structure} on $G^\cdot$ an element 
$\Delta \in Hom(G^\cdot,G^\cdot)$ of degree  $-1$ such that 
$d_H\Delta = [ , ]$ and $d_{CH}\Delta = 0$, i.e. such that
$$
\Delta(xy) - \Delta(x)y - (-1)^{\tx} x\Delta(y) = 
(-1)^{\tx}[x,y]
\eqno{(BV1)} 
$$
and 
$$
\Delta([x,y]) = [\Delta(x),y] + (-1)^{\tx - 1}[x,\Delta(y)],
\eqno{(qBV2)'}
$$
i.e. $\Delta$ is a derivation of the graded Lie algebra 
$G^{\cdot Lie} = G^\cdot [1]$. 

By 5.3, $(qBV2)'$  is  equivalent to  
$$
\Delta^2(xy) = \Delta^2(x)y + x\Delta^2(y)
\eqno{(qBV2)},
$$
i.e. $\Delta^2$ is a  est  derivation of the assocative algebra  $G^\cdot$. 

We see that $(qBV2)$ is weaker than the BV axiom, 
$\Delta^2 = 0$, cf. [K], \S 1.  

{\bf 5.5.} By \S 4, (a), the set  $qBV(G^\cdot)$ of qBV  structures 
on $G^\cdot$ is a torsor under the space 
$$
H^0(\text{Tot} C^{\cdot\cdot}_{HCH}(G^\cdot,G^\cdot)) = 
\text{Ker}(d^{00}_H) \cap \text{Ker}(d^{00}_{CH}) = 
$$
$$
\Omega^{1,\text{cl}}(G^\cdot) := \Ker\{d: \Omega^1(G^\cdot) \lra  \Omega^2(G^\cdot)\}.
$$

For $G^\cdot = \Lambda_A^\cdot(T)$ the $\Omega^{1,\text{fer}}(G^\cdot)$-torsor 
$qBV(G^\cdot)$ is identified with the  $\Omega^{1,\text{fer}}(T)$-torsor 
$BV(T)$. 

{\bf 5.6.} {\bf Remarks.} (a) Consider the bracket as an element  
$\omega \in C^{10}_{HCH}(G^\cdot,G^\cdot) = 
Hom(\Lambda^2(G^\cdot[1]),G^\cdot)[1]$,
$\omega(x,y) = [x,y]$. 

If $I \in C^{00}_{HCH}(G^\cdot,G^\cdot) = Hom(G^\cdot[1],G^\cdot[1])$ is the identity 
morphism then   $\omega = d_{CH}I$, whence 
$d_{CH}\omega = 0$. 

Computing the Hichschild differential,
$$
d_H\omega(x,y,z) = (-1)^{\tx\ty + \tx + 1}y[x,z] 
+ (-1)^{\tx + \ty}[x,yz] + 
(-1)^{\tx + \ty + 1}[x,y]z = 0,  
$$
i.e. $d_H\omega = 0$. Therefore 
$\omega\in \Omega^2(G^\cdot)^\sim$ ($\omega\notin \Omega^2(G^\cdot)$ since  $\tilde\omega =  - 1$), 
and $d_{DR}\omega = 0$. 

However, 
$$
d_HI(x,y) = (-1)^{\tx}xy, 
$$
so $\omega$ is not a coboundary in   $\Omega^2(G^\cdot)^\sim$. 

(b) If one defines  $m \in C^{01}(G^\cdot,G^\cdot)$ by $m = d_HI$, 
i.e. 
$$
m(x,y) = (-1)^{\tx}xy 
$$
then $d_Hm = 0$ and 
$$
d_{CH}m = d_{CH}d_H I = d_H d_{CH} I = d_H \omega = 0.
$$

{\bf 5.7.} It seems very probable that the contents of 
 [S2], Caput 2 (resp. of [S3]) can be generalized to the context of Gerstenhaber algebras,  i.e. 
one can define a notion of a {\it vertex structure} (resp. of a {\it membrane structure}) on a Gerstenhaber algebra. 
The set of these structures should form a $2$-torsor over a truncated 
 de Rham complex 
$$
\Omega^{[2,3>}(G^\cdot):\ \Omega^2(G^\cdot) \lra \Omega^{3,\text{cl}}(G^\cdot)
$$
(resp. over
$$
\Omega^{[3,5>}(G^\cdot):\ \Omega^3(G^\cdot) \lra\Omega^4(G^\cdot) \lra\Omega^{5,\text{cl}}(G^\cdot)).
$$

It would be interesting to understand the nature of the 
corresponding enveloping ''vertex Gerstenhaber algebras''. 
It is not immediate, due to the last term in Hochschild differential, 
which is not seen in the framework of vertex algebroids.

\bigskip\bigskip

\newpage

\centerline{\bf References}

\bigskip\bigskip

[BE] D.Bakry, M.\'Emery, Diffusions hypercontractives,  {\it S\'eminaire de Probabilit\'es} XIX 1983/84, 1985 - Springer.

[BGE] D.Bakry, I.Gentil, M.Ledoux, Analysis and geometry of  Markov 
diffusion operators, Springer, Grundlehren d. math. Wiss., 2014. 

[BD] A.Beilinson, V.Drinfeld, Chiral algebras, AMS Colloquium Publications, 
{\bf 51}, 2004. 

[CE] H.Cartan, S.Eilenberg, Homological algebra, Princeton University Press, 
1956.  

[GMS] V. Gorbounov, F. Malikov, V. Schechtman, Gerbes of  chiral 
differential operators. II. Vertex algebroids,{\it Inv. Math.}, {\bf 155}, 
605 - 680 (2004).  

[K] J.-L. Koszul, Crochet de Schouten - Nijenhuis et cohomologie, 
pp. 257 - 271 dans: Ast\'erisque, hors s\'erie, 1985, 
tome d\'edi\'e \`a E. Cartan, 1985. 

[S1] V. Schechtman, Remarks on formal deformations and Batalin - Vilkovisky
algebras, math.AG/9802006.   

[S2] V. Schechtman, Definitio nova algebroidis verticiani, pp. 443 - 494 
in: {\it Studies in Lie
theory: Antony Joseph's Festschrift}, Birkh\"auser, Progress in Math., {\bf 243}, 2006. 

[S3] V. Schechtman, Structures membranaires., pp. 599 - 643 in: {\it Algebraic geometry and Number theory, Vladimir Drinfeld's Festschrift}, 
Birkh\"auser, Progress in Math., {\bf 253}, 

[TT] D. Tamarkin, B. Tsygan, The ring of differential operators on forms 
in noncommutative calculus, pp. 105 - 132 in: {\it Graphs and Patterns in Mathematics
and Theoretical Physics,  Dennis Sullivan Festschrift}, Proc. Symp Pure Math., {\bf
73}, AMS, 2005. 

\bigskip\bigskip

Institut de Math\'ematiques de Toulouse, Universit\'e Paul Sabatier, 
31062 Toulouse, France

schechtman@math.ups-tlse.fr

\end{document}